\documentclass[11pt]{amsart}

\newtheorem{theorem}{Theorem}[section]

\theoremstyle{definition}
\newtheorem{remark}{Remark}[section]
\numberwithin{equation}{section}

\newcommand{\s}{{\mathcal S}}

\begin{document}

\title[Spatial Mean of Poincar\'{e} Cycle]
{On the Spatial Mean of the Poincar\'{e} Cycle}

\author{Luis B\'{a}ez-Duarte}

\date{\today}

\address{Departamento de Matem\'{a}ticas, 
Instituto Venezolano de Investigaciones
Cient\'{\i}ficas, Apartado 21827, Caracas 1020-A,
 Venezuela}

\email{lbaez@ccs.internet.ve}

\thanks{This article appeared in Spanish in the Bulletin of the Venezuelan Academy of Sciences in 1964 \cite{baez}, and is translated here into English by the author himself. Given the intervening 41 years, we deemed it convenient to make some minor improvements in style that change nothing essential in the content of the original paper.}


\begin{abstract}
This is an English translation of an article by the author appeared in Spanish in 1964. Let $X$ be a measure space and $T:X\rightarrow X$ a measurable transformation. For any measurable $E\subseteq X$ and $x\in E$, the possibly infinite return time is $n_E(x):=\inf\{n>0: T^n x\in E\}$. If $T$ is an ergodic tranformation of the probability space $X$, and $\mu(E)>0$, then a theorem of M. Kac states that 
$\int_E n_E d\mu=1$. We generalize this to any invertible measure preserving transformation $T$ on a finite measure space $X$, by proving independently, and nearly trivially that for any measurable $E\subseteq X$ one has $\int_E n_E d\mu=\mu(I_E)$, where $I_E$ is the smallest invariant set containing $E$. In particular this also provides a simpler proof of Poincar\'{e}'s recurrence theorem.
\end{abstract}
\maketitle

The purpose of this note is to give an independent and rather elementary proof of a slight generalization of Kac's well-known theorem \cite{kac} about the spatial mean of the Poincar\'{e} return time. The result implies Poincar\'{e}'s recurrence theorem as well. Let $(X,\s,\mu)$ be a finite measure space, and $T$ a measure preserving invertible transformation of $X$ onto itself, that is, 
$\mu(E)=\mu(T^{-1}E)=\mu(TE),\ \forall E\in\s$. For any $E\in\s$ define the associated function $n_E$ by
$$
n_E(x):=\inf \{n: n>0,\  T^n x\in E\},
$$
which may possibly be infinite; $n_E(x)$ is the \emph{first time of return of $x$ to $E$}, or \emph{the Poincar\'{e} cycle of $x$}. We denote by $I_E$ the smallest invariant set containing $E$. As usual $\chi_A$ is the characteristic function of the set $A$. 
\begin{theorem}
For every $E \in X$
\begin{equation}\label{spmean}
\int_E n_E  d\mu = \mu(I_E). 
\end{equation}
\end{theorem}

\begin{proof}
By definition of $n_E$ we have
\begin{eqnarray}\nonumber
\chi_{E^c}(Tx)\chi_{E^c}(T^2 x)\dots \chi_{E^c}(T^n x)
&=&
1, \ \ \ n<n_E(x),\\\nonumber
&=&
0, \ \ \ n\geq n_E(x),
\end{eqnarray}
which added over $n$ yields
\begin{eqnarray}\nonumber
n_E(x)
&=&
1+\sum_{n=1}^\infty 
\chi_{E^c}(Tx)\chi_{E^c}(T^2 x),\dots \chi_{E^c}(T^n x)\\\nonumber
&=&
1+\sum_{n=1}^\infty \chi_{\left(\bigcup_{\nu=1}^n T^{-\nu}E\right)^c}(x).
\end{eqnarray}
Now we integrate the above expression for $n_E$ to obtain
\begin{equation}\label{eq1}
\int_E n_E  d\mu = \mu(E)+\sum_{n=1}^\infty \mu\left(\left(\bigcup_{\nu=1}^n T^{-\nu}E\right)^c \cap E\right).
\end{equation}
But $T^n$ is invertible measure preserving so
$$
\mu\left(\left(\bigcup_{\nu=1}^n T^{-\nu}E\right)^c \cap E\right)=
\mu\left(\left(\bigcup_{\nu=0}^{n-1} T^{\nu}E\right)^c \cap T^n E\right),
$$
which substituted in (\ref{eq1}) yields
\begin{eqnarray}\nonumber
\int_E n_E  d\mu 
&=&
\mu(E)+\sum_{n=1}^\infty 
\mu\left(\left(\bigcup_{\nu=0}^{n-1} T^{\nu}E\right)^c \cap T^{n}E\right)
\\\label{eq2}
&=&
\mu\left(\bigcup_{\nu=0}^{\infty} T^{\nu}E\right),
\end{eqnarray}
where the last equality follows from the canonical disjoint decomposition of a countable union:
$$
\bigcup_{\nu=0}^{\infty} T^{\nu}E=
E \cup(E^c \cap TE)  \cup((E\cup TE)^c \cap T^2 E)\dots
$$ 
But it is obvious that $I_E=\bigcup_{n\geq0}T^n E$, which substituted in (\ref{eq2}) gives (\ref{spmean}).
\end{proof}
\begin{remark}
The fact that $\int_E n_E d\mu <\infty$ implies that $n_E(x)<\infty$ for a.e. $x$; this is Poincar\'{e}'s recurrence theorem in a slightly more general setting. For a probability space $X$ and an ergodic transformation $T$, if $\mu(E)>0$ then $I_E=X$, so we get Kac's theorem in \cite{kac}, that is $\int_E n_E d\mu = 1$.
\end{remark}

\bibliographystyle{amsplain}

\end{document}